\documentclass[12 pt]{article}
\usepackage{amssymb}
\usepackage[centertags]{amsmath}
\usepackage{amscd, amsthm}

\newtheorem{fact}{fact}[section]

\newtheorem{thm}[fact]{Theorem}
\newtheorem{lemma}[fact]{Lemma}

\newtheorem{corollario}[fact]{Corollary}
\newtheorem{defini}[fact]{Definition}
\newtheorem{osserva}[fact]{Remark}

\begin{document}
\title{From Schanuel's Conjecture to Shapiro's Conjecture}
\author{P. D'Aquino\footnote{Department of Mathematics, Seconda Universit\`{a} di Napoli, Via Vivaldi
43, 81100 Caserta, paola.daquino@unina2.it}, A. Macintyre
\footnote{School of Mathematical Sciences, Queen Mary, University
of London, Mile End Road, London E1 4NS, angus@dcs.qmul.ac.uk} and
G. Terzo\footnote{Department of Mathematics, Seconda
Universit\`{a} di Napoli, Via Vivaldi 43, 81100 Caserta,
 giuseppina.terzo@unina2.it}}

\maketitle

\begin{abstract}
In this paper we prove Shapiro's  $1958$ Conjecture on
exponential polynomials, assuming Schanuel's
Conjecture.
\end{abstract}

\section{Introduction}

We work with  exponential polynomial functions on $\mathbb C$ of the form
\begin{equation}\label{espon}
f(z) = \lambda_1 e^{\mu_1 z} + \ldots + \lambda_N e^{\mu_N z}.
\end{equation}
The set of such functions forms a ring  $\mathcal E$ under the usual addition and multiplication.
We normally refer to exponential polynomial functions simply as exponential polynomials.
In (\ref{espon}), we assume without loss of generality that the exponents $\mu$'s are distinct, and that the coefficients $\lambda$'s are nonzero, unless $f$ is the zero polynomial.

In  1974 during  the  Janos  Bolyai Society Colloquium on Number Theory, H. L. Montgomery
mentioned the following conjecture, which he attributed to H. S. Shapiro
\cite{shapiro conj}:\\

\noindent \textbf{Shapiro's Conjecture:} {If $f$ and $g$ are two exponential
polynomials in $\mathcal E$ with infinitely many common roots, then there exists
an exponential polynomial $h$ in $\mathcal E$ such that

$h$ is a common divisor of $f$ and $g$ in the  ring $\mathcal E$, and $h$ has infinitely many zeros in $\mathbb C$.


\bigskip

Montgomery pointed out, via an example given in \cite{vanderpoortentijdeman}, that the problem was not likely to yield easily to any classical approximation argument.

\bigskip
 It turns out that  Shapiro's Conjecture is naturally connected to Schanuel's Conjecture in Transcendence Theory.
\medskip

\noindent \textbf{Schanuel's Conjecture: } Let $\lambda_1, \ldots,
\lambda_n \in \mathbb C$. Then the transcendence degree of $\Bbb
Q(\lambda_1, \ldots, \lambda_n, e^{\lambda_1}, \ldots,
e^{\lambda_n})$ over $\Bbb Q$ is greater or equal than the linear
dimension of  $\lambda_1, \ldots, \lambda_n $ over $\Bbb Q$.

\medskip

Schanuel's Conjecture has played a crucial role in exponential
algebra (see \cite{angusfree}, \cite{GiusyEuler}), and in the
model theory of exponential fields (see \cite{MacWilkie},
\cite{Zilber}, \cite{Marker}).

\smallskip

In \cite{ritt} Ritt obtained a factorization theory for
exponential polynomials in $\mathcal E$. Subsequently, his ideas
have been developed, and his results have been extended to more
general exponential polynomials over $\mathbb C$, see
\cite{gourin}, \cite{macoll}, \cite{vanderpoorten},
\cite{vdpoorteneverest}. In \cite{paolagiusy}, these ideas have been put in the much broader context of general exponential polynomials (with any iteration of  exponentiation) over an
algebraically closed field of characteristic $0$ with an
exponentiation.

 In this paper we will study the  Shapiro Conjecture in a  context  more general than that of   the complex field. We will be working over an  algebraically closed field of characteristic $0$, with an exponential function, and  having an infinite cyclic group of periods, whose exponential  is surjective onto the multiplicative group.  The class of such fields includes the very important fields introduced by Zilber (see \cite{Zilber} for the basic notions).  The preceding assumptions play a minor role in our work on Shapiro's Conjecture. Of crucial importance is our further assumption, true for Zilber's fields, but unproved for the complex field, that we work with exponential fields satisfying Schanuel's Conjecture.

 One should note that in an exponential field satisfying the above assumptions (even without the surjectivity of the exponential onto the multiplicative group) the two element set of generators of the periods is first-order definable \cite{KMO}, the sine and cosine function are unambiguously defined, and the two element set consisting of the quotients of the period generators by   twice a square root of $-1$ is definable (the set does not depend on which root is chosen). In $\mathbb C$ this would define the set $\{ \pi, -\pi \}$. Finally, we can define the one element subset consisting of the element $x$ such that $sin(x/2) = 1$. In a general field satisfying our assumptions, we call this element $\pi$.


In section 2 we review the basic ideas of  Ritt's factorization theory for exponential polynomials.
His main theorem allows us to break the proof   of  Shapiro's conjecture into two cases. One case was already done by  van der Poorten and Tijdeman \cite{vanderpoortentijdeman} for simple polynomials (in Ritt's sense) over $\mathbb C$, without any use of Schanuel's Conjecture. In section 3, we modify that argument so as to apply to fields satisfying  all the assumptions given above with the exception of Schanuel's Conjecture.
\par
Section 4 explains recent work of Bombieri, Masser and Zannier \cite{bombieri} on anomalous subvarieties of powers of the multiplicative group.

 The main result of this paper is in section 5 where a positive solution to Shapiro's Conjecture is obtained for the
remaining case of irreducible exponential polynomials, assuming Schanuel's Conjecture,  using the work of Bombieri, Masser and Zannier, and work of
Evertse, Schlikewei and Schmidt on linear functions of elements of multiplicative groups of finite rank.

\par

We feel obliged to make a philosophical remark about the use of Schanuel's Conjecture to ``settle" a conjecture which emerged from complex analysis. A very distinguished number theorist has remarked that if one assumes Schanuel's Conjecture one can prove anything. The sense of this is clear if one restricts  ``anything" to refer to statements in transcendence theory. In that domain  the Conjecture is almost a machine, leading one mechanically to ``proofs" of  any plausible conjectures about algebraic relations between complex numbers and their exponentials. We suspect that there is also a common intuition that statements about common zeros of exponential poynomials should be related to statements about the transcendence theory of the exponential function. However, the original motivation for Shapiro's Conjecture clearly comes from reflection on distribution of zeros of individual exponential polynomials, and predates Schanuel's Conjecture. Moreover, our argument involves combinatorial considerations not previously connected to routine applications of Schanuel's Conjecture.


\section{Factorization Theory}\label{fattorizzazione}

We briefly review the main ideas in Ritt's factorization for  exponential polynomials
 in $\mathcal E$. Most of the theory adapts to the much more general context of
the ring of exponential polynomials over an  algebraically closed field of characteristic $0$ with an exponential function (see
\cite{paolagiusy}).

The fundamental idea  due to Ritt  was to transform  problems of factorization
 of  exponential polynomials to those of factorization of classical multivariate
polynomials in the extended category of polynomials in fractional powers of the variables. This brings in the notion of power
irreducible multivariate polynomial explained below.
\par In general, if we consider an irreducible polynomial $Q(x_1, \ldots, x_n)$
it can happen that for some positive integers $q_1, \ldots, q_{n}$ the
polynomial $Q(x_{1}^{q_1}, \ldots, x_n^{q_n})$ is reducible.

If there exists no sequence $q_1,
\ldots, q_{n}$ of positive integers such that $Q(x_{1}^{q_1}, \ldots, x_n^{q_n})$ is reducible we will refer to $Q$ as a {\it  power irreducible polynomial}.
\par We briefly review how to associate a classical polynomial in one or more variables
to an exponential polynomial in $\mathcal E$. \par We collect some basic definitions and results.
\medskip

\noindent {\bf Fact.}
The units in the ring $\mathcal E$ are the products of nonzero constants and $e^{\alpha z}$ for constant $\alpha \in \mathbb C$.

\begin{defini}
An element $f$ in $\mathcal E$ is irreducible, if there are no
non-units $g$ and $h$ in $\mathcal E$ such that $f = gh.$
\end{defini}


\begin{defini}\par Let $f = \sum_{i=1}^{N}\alpha_ie^{\mu_i z}$ be an exponential polynomial.
The {\it support} of $f$,  denoted by
$supp(f)$, is the $\mathbb Q$-space generated by
$\mu_{1}, \ldots, \mu_{N}$.\end{defini}

\begin{defini}
An exponential polynomial $f(z)$ of $\mathcal E$  is simple if
$\dim supp(f) = 1$.
\end{defini}

It is easily seen that, up to a unit, a simple exponential polynomial is a polynomial in
$e^{\mu z}$, for some $\mu \in \mathbb C$. An example of a simple
exponential  polynomial is $$g(z) =  \frac{e^{2\pi i z} -e^{-2\pi
i z}}{2i} = sin(2\pi z). $$

\begin{osserva}\label{simple} \rm{A simple exponential polynomial factorizes, up to units,  into a finite product of factors of the
form $1 - \alpha e^{\mu z},$ where $\alpha,\mu \in \Bbb C.$ This simply uses the fact that the complex field is algebraically closed. If
fractional powers of the variables are allowed then a simple
exponential polynomial may have infinitely many factors, e.g. $1 -
\alpha e^{\frac{\mu z}{k}},$ for each $k\in \mathbb N, k\not=0$.}
\end{osserva}

Let $f(z) = \lambda_1 e^{\mu_1 z} + \ldots + \lambda_N e^{\mu_N
z}$ where $\lambda_i$ and $\mu_i$ are complex numbers, and let
$\beta_1, \ldots \beta_r$ be a $\Bbb Z$-basis of the additive group generated by the $\mu_i$'s. Let
$Y_j = e^{\beta_j z},$ with $j = 1, \ldots, r.$ If each $\mu_i$ is
expressed in terms of the $\beta_j$'s we have that $f$ is
transformed into a classical Laurent  polynomial $Q$ over $\Bbb C$ in the
variables $Y_1, \ldots, Y_r$. The best way to think of $Q$  is as a   function
on the product of $r$ copies of the multiplicative group variety. Remember that the
 $Y$'s   are exponentials and so take value in the multiplicative group. More prosaically, one can write $Q$ as a
product of a polynomial in the $Y$'s and a quotient of monomials in the $Y$'s.
\par Clearly, any factorization of $f$ determines a
factorization of $Q(Y_1, \ldots, Y_r).$ Ritt saw the relevance, in
terms of factorization theory, of understanding the ways in which
an irreducible polynomial $Q(Y_1, \ldots, Y_r)$ can become
reducible once  the variables are replaced by their powers. It is a fundamental problem to determine the set of integer  $r$-tuples $q_1, \ldots ,q_r$ for which the reducibility
occurs. Ritt gave a uniform bound for the number of
irreducible factors of $Q(Y_1^{q_1}, \ldots, Y_r^{q_r})$, depending
only on the
degree of $Q.$\\

\par For the factorization theorem of  Ritt the following lemma
is crucial.

\begin{lemma}
\label{divisibilitycondition}
Let $f(z)=\sum_{i=1}^{N}\alpha_ie^{\mu_i z}$ and $g(z) = \sum_{j=1}^{M}l_je^{m_j z}$ be non-zero exponential polynomials. If $f$  is
divisible by $g$ then  $supp(ag)$ is
contained in $supp(bf)$, for some units $a$ and $b$,  i.e. every element of $supp(ag)$ is a linear
combination of elements of  $supp(bf)$ with rational coefficients.
\end{lemma}

Note that if $f$ is a simple polynomial and $g$ divides $f$ then $g$ is also simple.  The factorization theorem that we need is the following (see \cite{ritt}, \cite{gourin} and \cite{macoll}).

\begin{thm}
\label{Ritt} Let $f(z) = \lambda_1 e^{\mu_1 z} + ... + \lambda_N
e^{\mu_N z}$, where $\lambda_i, \mu_i \in \mathbb C$. Then $f$ can
be written uniquely  up to order and multiplication by units as $$f(z) = S_1 \cdot \ldots \cdot S_k \cdot
I_1 \cdot \ldots \cdot I_m
$$ where $S_j$ are simple polynomials with $supp(S_{j_1})\not=
supp(S_{j_2})$ for $j_1\not= j_2$, and $I_h$ are irreducible
polynomials in $\mathcal E$.
\end{thm}

We observe that the proof has nothing to do with analytic functions, and works over any characteristic
 $0$ exponential field which is  an algebraically closed  field (see \cite{vdpoorteneverest}, \cite{paolagiusy}). This is the context where we will be working.
 \par
 Since a common zero of two products is a common zero of two factors,
Theorem \ref{Ritt} trivially implies that only two cases of the  Shapiro
Conjecture have to be considered.

\smallskip

\noindent {\bf Case 1.} At least one of the exponential polynomials
$f$ and $g$ is simple.

\smallskip
\noindent {\bf Case 2.} Both of the exponential
polynomials $f$ and $g$  are  irreducible.

\section{Shapiro Conjecture: Case 1}

The case when either $f$ or $g$ is simple has been proved unconditionally
 by van der Poorten and Tijdeman for the complex field, see
\cite{vanderpoortentijdeman}.

Their proof  uses various results from  Ritt divisibility theory
in \cite{ritt1} and  a variant of the usual $p$-adic  argument
from the proof of the Skolem-Mahler-Lech  Theorem on recurrence
sequences with infinitely many vanishing terms. Ritt's result most
specific to the complex field says  that if    $f/g$ is an entire
function, where $f$ and $g$ are exponential polynomials, then $f$
divides $g$. The proof ultimately relies on a fundamental result
of Tamarkin, Polya and Schwengler on the distribution of zeros for
exponential polynomials as in (\ref{espon}).  We observe that it
is not obvious what interpretation to give this result in more
general exponential fields, and for that reason we have sought and
found a proof that avoids this result of Ritt. We do not, however,
avoid appeal to the Skolem-Mahler-Lech Theorem. The latter
theorem, as used in \cite{vanderpoortentijdeman} on Page $62$, in
a formulation for exponential functions, is:

\begin{thm}\label{SML}{\rm{(Skolem, Mahler, Lech)}} If  $f(z)$ is a function
as in (\ref{espon}) which vanishes for infinitely many
integers $z$ then there exists an integer $\Delta$ and positive residues $d_{1}, \ldots,
d_{l}$ modulo $\Delta$, such that $f(z)$ vanishes for all integers $z \equiv
d_{i}(\rm{mod }\mbox{ } \Delta),$ $i = 1, \ldots, l,$ and $f(z)$
vanishes only finitely often on other integers.
\end{thm}
Inspection of the proof (by a suitable $p$-adic embedding) shows
that it works for all exponential fields of characteristic $0$.
\par

We extend  the van der Poorten - Tijdeman result to the more
general setting of an exponential algebraically closed field $K$
of characteristic $0$, with standard periods and  exponential map
surjective to the multiplicative group, making no use of analytic
methods.  We need the following lemma.

\begin{lemma}
\label{sine} Let $h(z)=\lambda_1 e^{\mu_1 z} + \ldots + \lambda_N
e^{\mu_N z},$ where $\lambda_j, \mu_j\in K$. If $h$ vanishes at
all integers then $sin(\pi z)$ divides $h$.
\end{lemma}

\noindent
\emph{Proof:} We proceed by induction on the length $N$ of $h$. If $N=2$ the proof is a trivial direct computation.
\par
Let $N>2$, and consider the first $N$ positive solutions $1, \ldots, N$. The following identities hold \\
$$\left\{ \begin{array}{lll}\lambda_{1}e^{\mu_{1}} + \lambda_{2}e^{\mu_{2}} + \ldots
\ldots +
\lambda_{N}e^{\mu_{N}} & = 0 &\\
& & \\ \lambda_{1}(e^{\mu_{1}})^2 + \lambda_{2}(e^{\mu_{2}})^2 +
\ldots \ldots +
\lambda_{N}(e^{\mu_{N}})^2 & = 0 &\\
\cdots & &\\
\cdots & &\\
\lambda_{1}(e^{\mu_{1}})^N + \lambda_{2}(e^{\mu_{2}})^N + \ldots +
\lambda_{N}(e^{\mu_{N}})^N & = 0 &
\end{array} \right.$$\\

Let $\delta_{1} = e^{\mu_{1}}, \cdots, \delta_{N} =
e^{\mu_{N}},$ so by substitution we can rewrite the identities in
matrix notation  as follows\\
$$\left( \begin{array}{llcl} \delta_{1} & \delta_{2} &\ldots &\delta_{N}\\
\delta_{1}^{2}& \delta_{2}^{2} & \ldots & \delta_{N}^{2}\\
\vdots & \vdots & \ddots & \vdots\\
\delta_{1}^{N} & \delta_{2}^{N} & \ldots & \delta_{N}^{N}
\end{array}\right)\left(
\begin{array}{l} \lambda_1\\ \lambda_2\\ \vdots\\ \lambda_N \\\end{array}\right)
=\left(
\begin{array}{l} 0\\ 0\\ \vdots\\ 0 \\\end{array}\right).$$

Because of the existence of a non trivial solution of the system the determinant of the matrix vanishes,

$$\left| \begin{array}{llcl} \delta_{1} & \delta_{2} &\ldots &\delta_{N}\\
\delta_{1}^{2}& \delta_{2}^{2} & \ldots & \delta_{N}^{2}\\
\vdots & \vdots & \ddots & \vdots\\
\delta_{1}^{N} & \delta_{2}^{N} & \ldots & \delta_{N}^{N}
\end{array}\right|
 = 0,$$
 that is\\

$$\delta_{1} \cdot \delta_{2}\cdot \ldots \cdot \delta_{N} \cdot \left| \begin{array}{llcl} 1 & 1 & \cdots & 1\\
\delta_{1}& \delta_{2} & \ldots & \delta_{N}\\
\vdots & \vdots & \ddots & \vdots\\
\delta_{1}^{N-1} & \delta_{2}^{N-1} & \ldots & \delta_{N}^{N-1}
\end{array}\right|
 = 0.$$

This is a  Vandermonde determinant, so:\\

$$(\delta_{1} \cdot \delta_{2} \cdot \ldots \cdot \delta_{N}) \cdot
\prod_{1\leq i < \ell \leq N}(\delta_{i} - \delta_{\ell}) = 0.$$
So,  $\delta_{i}=\delta_{\ell}$ for some $i \neq \ell,$ i.e.
$e^{\mu_{i}} = e^{\mu_{\ell}}$ for some $i \neq \ell,$ and without loss
of generality we can assume $e^{\mu_1 }=e^{\mu_2}$. So,  $e^{\mu_1n }=e^{\mu_2n}$
for each $n\in \mathbb Z$. The polynomial $$(\lambda_1+ \lambda_2)e^{\mu_1 z} + \sum_{j\geq 3} \lambda_j e^{\mu_j z}$$
also vanishes on all integers, and since it has
length strictly less than $N$ it is divisible by $sin(\pi z)$.
Note that $$h(z)=(\lambda_1+ \lambda_2)e^{\mu_1 z} + \sum_{j\geq 3}+ \lambda_j e^{\mu_j z}+ \lambda_2(e^{\mu_1z }-e^{\mu_2z}),$$
and so all integers are roots of $e^{\mu_2z }-e^{\mu_1z}$.
This implies that $\mu_2-\mu_1\in 2\pi i\mathbb Z$,
hence $$e^{\mu_1z }(e^{2\pi ijz }-1)=e^{\mu_1z }e^{\pi ijz }(e^{\pi ijz }-e^{-\pi ijz })$$
for some $j\in \mathbb Z$, and clearly $sin(\pi z)$ divides $e^{\pi ijz }-e^{-\pi ijz }.$ \hfill $\Box$ \\

The Shapiro Conjecture for the case when one of the polynomials is
simple follows from the following theorem which implies that if
one of the two polynomial is simple so is  the other one.
\begin{thm}
Let $f$ be a simple exponential polynomial, and let $g$ be
an arbitrary exponential polynomial such that $f$ and $g$ have infinitely many common roots.
Then there exists an exponential polynomial which divides both $f$ and $g$.
\end{thm}

\noindent {\it Proof:} If $f$ is simple then up to a constant, $f$ is of the form,
$f=\prod (1-ae^{\alpha z})$, where $a, \alpha \in K$. If $f$ and $g$ have infinitely
 common zeros then $g$ has infinitely common zeros with one factor of $f$, say $1-ae^{\alpha z}$.
So $g$ has infinitely many zeros of the form $z=(2k\pi i-\log
a)/\alpha$ with $k\in \mathbb Z$, and for a fixed value of $\log
a$. If $g^*(z)=g((2\pi iz-\log a)/\alpha) $ then $g^*$ has
infinitely many zeros in $\mathbb Z$. By Theorem \ref{SML},
$g^*(z)$ vanishes on  the set $M=\{ d_0+j\Delta :j\in \mathbb
Z\}$, for some $\Delta $ and $d_0$ in $\mathbb Z$, and $0\leq
d_0<\Delta$. If $h(z)=g^*(d_0+z\Delta)$ then $h$ vanishes on
$\mathbb Z$, and Lemma \ref{sine} implies that $h$ is divisible by
$sin(\pi z)$. This is a contradiction if $h$ is irreducible, which
is the case when $g$ is irreducible. This forces $g$ to be simple
(up to a unit), e.g. $g(z)=1-be^{\beta z}$ for some $b,\beta \in
K$.  So, without loss of generality we can consider the system
\begin{equation}\label{sistemaridotto}
\left\{ \begin{array}{lll} f(z) = 1-ae^{\alpha z} & = 0 &\\
 g(z) = 1-be^{\beta z}& = 0 &
\end{array} \right.\\
\end{equation}
where $a,b,\alpha ,\beta \in K$, with infinitely many common
zeros. The roots of $f$ are of the form $z=\frac{1}{\alpha}(-\log
a+2k\pi i)$, $k\in \mathbb Z$. It follows that $g$ vanishes on
$z=\frac{1}{\alpha}(-\log a+2t\pi i)$ for infinitely many $t$ in
$\mathbb Z$. We argue now as before, using Theorem \ref{SML}, to
conclude that  $f$ and $g$ vanish on $$ \frac{1}{\alpha}(-\log a
+2(d+\Delta j)\pi i),$$ where  $d, \Delta $ are integers $d<
\Delta $, and for all $j\in \mathbb Z$. Via the change of variable
$T(z)= \frac{1}{\alpha}(-\log a +2(d+\Delta z)\pi i)$ the
exponential polynomials $f(T(z))$ and $g(T(z))$ both vanish on
$\mathbb Z$, and by Lemma \ref{sine} they are both divisible by
$sin(\pi z)$. Thus $f(w)$ and $g(w)$ are both divisible by $sin
(\pi T^{-1}(w))=sin (\frac{1}{2\Delta i}(\alpha w+\log a-2d\pi
i))$ which is a simple polynomial.  \hfill $\Box$

\section{Group varieties associated to exponential polynomials}

We now adapt, to the system $f=g=0$ the procedure of Ritt, thereby
converting the system to one defined by two conventional
polynomials, defining a subvariety of a power of the
multiplicative group.

\par
We work over an algebraically  closed  characteristic $0$
exponential field $K$ with standard periods,  with the exponential
surjective onto the multiplicative group, and satisfying
Schanuel's Conjecture ($SC$).
\\
Consider a  system with no restriction on $f$ and $g$

\begin{equation}\label{siste}
\left\{ \begin{array}{lll} f(z) = \lambda_1e^{\mu_1 z} + \ldots + \lambda_Ne^{\mu_N z} & = 0 &\\
 g(z) = l_1e^{m_1 z} + \ldots + l_Me^{m_M z} & = 0 &
\end{array} \right.\\
\end{equation}
where $\lambda_{i},\mu_i,l_j,m_j \in K.$\\

Let $D$ be the linear dimension of $supp(f)\cup supp(g)$, and
$b_{1}, \ldots, b_{D}$ a $\Bbb Z$-basis of the group generated by
$ \overline{\mu}, \overline{m} $.  We introduce new variables $Y_1
= e^{b_1z}, \ldots, Y_D = e^{b_Dz},$ and as in Section
\ref{fattorizzazione} we associate the Laurent polynomials
$$F(Y_1, \ldots, Y_D), G(Y_1, \ldots, Y_D) \in \Bbb
Q(\overline{\lambda}, \overline{l})[Y_1, \ldots, Y_D]$$ to $f(z)$
and $g(z),$ respectively. As far as zeros from the multiplicative
group are concerned, one may replace $F$ and $G$ by ordinary
polynomials got by multiplying them by monomials. Note that $F$
and $G$ are polynomials over $\mathbb
Q(\overline{\lambda},\overline{l})$. Let $L$ be the algebraic
closure of this field. Obviously $L$ has finite transcendence
degree, a fact which will be crucial later.

Clearly, if $s$ is a common zero of $f$ and $g$ then $(e^{b_1s},
\ldots, e^{b_Ds})$ is a common zero of $F$ and $G$ in the $D$th
power of the multiplicative group. The study of the set of
solutions of system (\ref{siste}) will be reduced to studying the
solutions of system
\begin{equation}
\label{associatedsystem}
\left\{ \begin{array}{lll} F(Y_1, \ldots, Y_D) & =& 0\\
 G(Y_1, \ldots, Y_D) & = & 0
\end{array} \right.\\
\end{equation}

\begin{osserva}
\label{dimensionvariety}
\rm{Let $V(F)$ and $ V(G)$ be the subvarieties in the $D$th power of the multiplicative group
$G^{D}_m$,  associated to $F$ and
$G,$ respectively.
 If $f$ and $g$ are irreducible then $F$ and $G$
are power irreducible. In this case $\dim V(F) = \dim V(G) = D -
1.$ If we assume that $f$ and $g$ are distinct irreducibles (i.e.
neither is a unit times the other) then $F$ and $G$ are power
irreducibles,   with neither a scalar multiple of the other. It
follows that the algebraic set defined by $F=G=0$ has dimension no
more than $D-2$. This is crucial in what follows.}

\end{osserva}
Recall that an algebraic  subgroup in the group variety $G^D_m$ is
given by a  finite set of conditions  each of the form
$$Y_1^{a_1}\ldots Y_D^{a_D}=1$$ where $a_1,\ldots ,a_D\in \mathbb
Z$.  We will refer to $(a_1,\ldots ,a_D)\in \mathbb Z^D$ as the
exponent vector. For such a variety, the dimension is $D-h$ where
$h$ is the rank of the subgroup of $\mathbb Z^D$ generated by the
exponent vectors. A translate or coset  of a subgroup is obtained
by replacing $1$  by  other constants in the finite set of
conditions. A torus is a connected algebraic subgroup.

\par
The algebraic set $C$ defined by (\ref{associatedsystem}) may be a
reducible subvariety of the algebraic group $G_m^D$ over $L$. As
remarked above, its dimension is at most $D-2$ if $F$ and $G$ are
distinct irreducible over $L$.

Later, in the proof of the Shapiro Conjecture, we will work on a suitable irreducible component of $C$.
\par Note that if $f$ and $g$ have infinitely many common zeros, and $f$ is irreducible,
the algebraic set $C$ above cannot be contained in any coset of any proper algebraic subgroup of $G_m^D$.
For otherwise, let
\begin{equation}
\label{coset}
Y_1^{a_1}\ldots Y_D^{a_D}=\theta
\end{equation}
be one of the equations defining the coset. This corresponds to a
simple polynomial in Ritt's sense which has infinitely common
zeros with $f$. By van der Poorten and Tijdeman result and Lemma
\ref{divisibilitycondition} we have a contradiction since $f$ is
not simple.

\medskip


We now review the basic concepts concerning the notion of
anomalous subvariety, as used in \cite{bombieri} by Bombieri,
Masser and Zannier.  We will not give the full details of the
analysis obtained by Bombieri, Masser and Zannier but we will
describe those properties of anomalous varieties which we will
need  in the proof of our main result. Their discussion  is first
done over the complexes, but they observe that it works over any
algebraically closed field of characteristic 0, and we use this
fact. For us the case of the $L$ introduced earlier is crucial because of its finite transcendence degree. We will follow \cite{bombieri} for the notion of a subvariety of the algebraic group $G_m^n$, and when necessary we will specify if the variety is irreducible.
\par

Let $V$ be an irreducible subvariety of $G_m^n$.

\begin{defini}
An irreducible subvariety $W$ of $V$ is anomalous in $V$
if $W$ is contained in a coset of  an algebraic subgroup $\Gamma$ of $G^{n}_m$ with

$$\dim W > max\{0, \dim V - {\rm codim }\Gamma \}$$
\end{defini}

Note that this definition has the same meaning in any algebraically closed field over which $V$ is defined.

\begin{defini} An anomalous subvariety of $V$ is maximal if it is
not contained in a strictly larger anomalous subvariety of $V.$
\end{defini}

\begin{thm}
\label{BMZ} Let $V$ be an irreducible variety in
$G^{n}_m$ of positive dimension defined over $\Bbb C.$ There is a
finite collection $\Phi_{V}$ of  proper tori $H$ such that
$1\leq n-\dim H\leq \dim V$ and every maximal anomalous subvariety
$W$ of $V$ is a component of the intersection of $V$ with  a coset $ H\theta$ for some $H \in
\Phi_{V}$ and  $\theta \in G^{n}_m$.
\end{thm}

For the proof see \cite{bombieri}. Note that this result is true
(as is stated in \cite{bombieri}) when $\Bbb C$ is replaced by any
algebraically closed field $K$ of characteristic $0$, in the sense
that the cosets involved, for $W$ defined over $K$, are also
defined over $K$.

Theorem \ref{BMZ} implies, since every anomalous subvariety is contained in a maximal one, that there is a finite number of
subgroups of codimension $1$, such that any anomalous subvariety is included in a coset of one of them.
\par


\section{The Full Shapiro Conjecture}
\par We concentrate now on  Case 2  of Shapiro's conjecture. In this case the conjecture has the following formulation:
{\it If $f$ and $g$ are distinct irreducible exponential polynomials 
then  $f$ and $g$ have at most finitely many common zeros.} 
\par
We
will prove the following equivalent version (see
\cite{vanderpoortentijdeman}): {\it Let $f$ and $g$ be   exponential polynomials, and assume $f$ is irreducible.  If $f$ and $g$ have
infinitely many common zeros then $f$ divides $g$.}

In the following unless otherwise specified the
linear dimension and the transcendence degree of a tuple will always be over $\Bbb Q.$

\smallskip

Let $D=l.d.(supp(f)\cup supp(g))$, and let $b_{1}, \ldots, b_{D}$
a $\Bbb Z$-basis of the group generated by $ \overline{\mu},
\overline{m}$. We will denote the transcendence degree of
$\overline \lambda, \overline l$ by $\delta_1$,  and the
transcendence degree of $\overline \mu, \overline m$ by
$\delta_2$, i.e. $\delta_1 = t.d.(\overline \lambda, \overline
l),$ and $\delta_2 = t.d.(\overline \mu, \overline m).$ We denote
by $\overline b$ the sequence $(b_1, \ldots, b_D)$ and by $B$ the
set $\{b_1, \ldots, b_D\}.$
\par

Assume that $f$ and $g$  have infinitely many common zeros. Let
$S$ be an infinite set of  nonzero common solutions. We will
``thin" this set inductively to infinite subsets using arguments
of Schanuel type, and work of Bombieri, Masser and Zannier on
anomalous intersections, to reach an infinite $S$ such that the
$\Bbb Q$-space generated by $S$ is finite dimensional. We will
then get a contradiction from using, inter alia, work of Evertse,
Schlickewei and Schmidt on linear functions of elements of finite rank groups.

\par We begin with some simple bounds on Schanuel data. For any $s
\in S$ let $\overline{b}s$ stand for the sequence $(b_1s, \ldots,
b_Ds)$
and  $e^{\overline bs}$
stand for the sequence $(e^{b_1s}, \ldots, e^{b_Ds})$. In terms of
the set, for any $s \in S$ we denote by $Bs = \{b_1s, \ldots,
b_Ds\},$ and by $e^{Bs} = \{e^{b_1s}, \ldots, e^{b_Ds}\}.$ For any
subset $T$ of the set of solutions $S,$ $$BT = \bigcup_{s \in T}
Bs,$$ and $$e^{BT} = \{e^{b_is} : 1 \leq i \leq D, b_i\in B, s \in T\}.$$
\par For any finite subset $T$ of $S$, let $D(T)$ be the linear dimension
of the space spanned by $BT.$
Notice that $D(T)=D$ if $T$ is a singleton, since $0\not\in S$.
Moreover, $D(T)\leq D|T|$, where $|T|$ denotes the cardinality of
$T$. We show now that there is an upper bound to the cardinality
of $T$ for which the equation $D(T)=D|T|$ holds.


\begin{lemma}
\label{maximum} (SC) For any finite subset $T$ of $S$ with
$D(T)=D|T|$ we have that $|T| \leq \delta_1 + \delta_2.$ 
\end{lemma}

\noindent \emph{Proof:} Enumerate the set $T$ as $s_1,\ldots ,s_k,$ of
elements of $S$. By previous observations, upper bounds on the respective 
transcendence degrees of the sets $e^{BT}$ and $BT$
are $$t.d.(e^{BT})\leq k(D-2) + \delta_1,$$(because of the
dimension estimate on $F=G=0$ given in Remark
\ref{dimensionvariety}) and
$$t.d.(BT) \leq \delta_2+k.$$
By Schanuel's Conjecture we have
$$t.d.(BT, e^{BT}) \geq D(T),$$
and this implies
\begin{equation}\label{schanuel}
D(T)\leq kD-k+\delta_1+\delta_2.
\end{equation}

If $D(T) = kD$, inequality (\ref{schanuel}) implies that

\begin{equation}\label{k}
\delta_1 + \delta_2 \geq k,
\end{equation}
for all $k \in \Bbb N,$
proving the result since $\delta_1$ and $\delta_2$ are fixed
and depend only on the coefficients of the polynomials $f$ and
$g.$
\qed

\begin{osserva}
\label{remarkmaximum} \rm{Let $k_0$ be the maximum cardinality of
a $T$ for which the equation $D(T) = D|T|$ holds. Let $S_0$ be such a $T$. If we extend $S_0$ to a set $S_1$,
by adding $k_1$ distinct elements, then we clearly have the
following estimates:

$$D(S_0) \leq D(S_1) \leq \delta_1 + \delta _2 + k_1(D-1).$$}


\end{osserva}


Lemma \ref{maximum} has a fundamental consequence on the
transcendence degree of the set $BS$ which will be crucial in the
following.

\begin{lemma}
\label{finitetd} (SC) The transcendence degree  of
$BS$ over $\Bbb Q$ is less or equal than $\delta_1 + 2\delta_2.$

\end{lemma}

\noindent 
\emph{Proof:} Fix any $s\in S-S_0$. Then by maximality of $S_0$
for the equation $D(T) = D|T|,$ we have a nontrivial linear
function $\Lambda$ over $\Bbb Q,$ such that $\Lambda(\overline{b}
s)$ belongs to the $\mathbb Q$-vector space generated by the
$BS_0.$ We note that $\overline b$ is linearly independent over
$\Bbb Q$ and $\Lambda$ is linear, so we get that $$s =
\Lambda(\overline b)^{-1} \cdot a$$ where $a$ is in the $\mathbb
Q$-vector space generated by the $BS_0.$ Let $F$ be the field
generated by $BS_0 \cup B.$



The transcendence degree of $F$ is clearly finite, and the
following inequalities hold
$$t.d._{\mathbb Q} (F)\leq t.d.(B)+t.d._{\mathbb Q(B)}
(BS_0)\leq \delta_2+k_0\leq
\delta_1+\delta_2+\delta_2=\delta_1+2\delta_2.$$
\qed




\bigskip

The following result will be crucial  for completing the proof of Shapiro's Conjecture.

\

\noindent
{\bf Main Lemma.}
{\it (SC) For some infinite subset $S'$ of S the $\Bbb Q$-vector space generated by $S'$ is finite
dimensional.}

\

\noindent  \emph{Proof:}
Consider the subvariety  $C$ of  $G_m^D$  defined by
\begin{equation}\label{sistema1}
\left\{ \begin{array}{lll} F(Y_1, \ldots, Y_D) & =& 0\\
 G(Y_1, \ldots, Y_D) & = & 0
\end{array} \right.\\
\end{equation}
over $L=\mathbb Q(\overline{\lambda}, \overline{l})^{alg}.$ This may be
a reducible subvariety of the algebraic group $G_m^D$, so we work
now with a fixed irreducible component  $V$  of $C$ containing
solutions of the form $(e^{b_1s}, \ldots, e^{b_Ds}),$ for
infinitely many  $s \in S.$ 

An upper bound on the dimension of $V$ over $L$  is $D-2$, and  so
$D-2+\delta_1$ is the corresponding upper bound  over $\mathbb Q$ (see Remark \ref{dimensionvariety}).

\par

We now thin $S$ to an infinite subset  $S^{'}$ such that for $s \in S^{'}$,  the $D$-tuple $e^{\overline b s}$ is a point of $V$. This might force to throw out part of the original $S_0$ but this is irrelevant for the estimates on the linear dimension of $S^{'}$.


\par Fix a finite sequence  $\overline s = (s_1, \ldots s_k)$ of distinct elements of
$S,$ of length $k,$ and let $T$ be the set of entries $\overline
s.$ The $\mathbb Q$-linear relations among ${\overline{b}s_1},
\ldots ,{\overline{b}s_{k}}$ can be converted into $\mathbb
Z$-linear ones, and these naturally induce multiplicative
relations of group type among the corresponding exponentials
$e^{\overline{b}s_1}, \ldots ,e^{\overline{b}s_{k}}$. Thus we
determine an algebraic subgroup $\Gamma_{k}$ of $G_m^{Dk}$ on
which $e^{\overline{b}s_1}, \ldots ,e^{\overline{b}s_{k}}$ lie.
Clearly, the codimension of $\Gamma_{k}$ is $Dk-D(T)$, and
dimension of $\Gamma_{k}$ over $\mathbb Q$ is $D(T).$

Let $V^{k}$ be the product variety in the multiplicative group
$G_m^{Dk}$. The $Dk$-tuple
\begin{equation}
\label{tuple} (e^{\overline{b}s_1}, \ldots ,e^{\overline{b}s_{k}})
\end{equation}
lies on it, and this is true for any choice of $k$ solutions $s_1,
\ldots ,s_k$. An upper bound for the transcendence degree of any
tuple as in (\ref{tuple}) over $L$ is $k(D-2)$, and $k(D-2)+\delta
_1$ is a corresponding upper bound over $\mathbb Q$.

\medskip
The $Dk$-tuple $$(e^{\overline{b}s_1},
\ldots,e^{\overline{b}s_{k}})$$ belongs to the intersection of
$V^{k}$  and $\Gamma_{k}$, which might be reducible, and we will
work with the variety $W_{\overline s}$ of the point
$(e^{\overline{b}s_1}, \ldots, e^{\overline{b}s_{k}})$ over $L.$

\

\noindent
{\bf Claim 1.} For $k>\delta_1 +\delta_2$ the variety $W_{\overline s}$ is either anomalous or of dimension $0$ over $L$. 
\par

Suppose $dim(W_{\overline s})\leq dim (V^{Dk})-codim
(\Gamma_{k})$, i.e. $dim(W_{\overline s}) \leq k(D-2)-
(kD-D(T))+\delta_1$.
Again Schanuel's Conjecture implies
$$D(T) \leq k(D-2)- (kD-D(T))+\delta_1+
\delta_1+2\delta_2 ,$$ and so $$2k\leq 2\delta_1+2\delta_2.$$
Hence the claim is proved.

\

\par We want to get results not sensitive to any particular enumeration. Now suppose we rearrange the sequence $\overline s$ to $\overline {s^{*}}$.
The set $T$ does not change. It is easy to see that we still get points on
$V^{k}$, and dimension of $\Gamma_k$ does not change. What may change is
$W_{\overline {s^{*}}}.$  But consider the automorphisms (of affine $Dk$-space, of $V^k$, and of $G_m^{Dk}$) got by simply permuting the natural $D$-blocks. These transform the $W_{\overline s}$ to the $W_{\overline {s^{*}}}$, and  one sees easily that
$W_{\overline s}$ has dimension $0$ if and only if  $W_{\overline
{s^{*}}}$ has, and that $W_{\overline s}$ is anomalous if and only
if $W_{\overline {s^{*}}}$ is. So the claim implies that for every $k$, if
$k > \delta_1+\delta_2$ then either each $W_{\overline
{s^{*}}}$ has dimension 0, or each $W_{\overline {s^{*}}}$ is
anomalous.


\

\noindent
{\bf Claim 2.} If $\dim W_{\overline s} = 0$ then $D(T) \leq 2\delta_1 +
2\delta_2.$

Suppose $\dim W_{\overline s} = 0.$ Hence the coordinates
of all elements of $W_{\overline s}$ are algebraic over $L$, which
implies that
$$t.d.(e^{\overline{b}s_1}, \ldots, e^{\overline{b}s_{k}}) \leq \delta_1.$$
From Lemma \ref{finitetd} it follows that 
$$t.d.(\overline{b}s_1, \ldots, \overline{b}s_{k}) \leq \delta_1 + 2\delta_2,$$
and Schanuel's Conjecture implies  $$D(T) \leq 2\delta_1 +
2\delta_2.$$


\par This now gives that  for any  $k$ element subset $T$ of $S,$ if $$D(T) > 2\delta_1 + 2\delta_2$$ then $W_{\overline s}$
is anomalous, for any enumeration $\overline s$ of $T.$





\medskip

\medskip

We consider now a countably infinite subset of  $S$ enumerated as
$s_1, s_2,\ldots,$ which we will continue to call $S$.  Define $S_k$ as the set $\{s_1, \ldots,
s_k\}.$ Let $W_k$ be one of the $W_{\overline s}$ for a sequence
$\overline s$ enumerating $S_k.$ If infinitely many $W_k$ are of
dimension 0, then the set $\{ e^{b_j s}:  b_j\in B, s\in S\}$ is contained in $L$, and so by Schanuel conjecture and
the preceding calculations, $D(S_k) \leq 2(\delta_1 + \delta_2)$ for infinitely many $k$'s.
So $S$ spans a finite dimensional space over $\Bbb Q$, which is
the required conclusion. Thus, there is a $k_1$ such that for $k$
at least $k_1$ no $W_k$ has dimension 0. Thus by Claim 1,
all $W_k$ are anomalous. Since $W_k$ was chosen for an arbitrary
enumeration of $S_k$, we conclude that each $W_{\overline {s^*}}$
is anomalous, for any enumeration $\overline {s^*}$ of $S_k.$\\

We will make use of the Bombieri, Masser and Zannier
results. Though the $W_k$ are defined relative to an enumeration,
and would change if the enumeration did, there are some basic
results independent of the enumeration, and these will be needed
in the remaining stages of the proof.

Let $k_2\in \mathbb N$ be the least integer  $k$ such that for any
$k_2+1$ elements of $S$, $\eta_1, \ldots , \eta_{k_2+1}$,  the
variety $W$ of  the  $k_2+1$-tuple
$e^{\overline{b}\eta_1}, \ldots , e^{\overline{b}\eta_{k_2+1}}$ is
anomalous in $V^{k_2 + 1}.$
From \cite{bombieri} it follows that there is a finite collection
$\Phi_{V^{k_2+1}}$ of proper tori $H_1, \ldots , H_t$ of
$G_m^{(k_2+1)D}$ such that each maximal anomalous subvariety of
$V^{k_2+1}$ is a component of the intersection of $V^{k_2+1}$ with a coset of one of the $H$'s.
\par
We use a much less precise version for general anomalous
subvarieties. This version follows from the very precise Structure
Theorem of Bombieri, Masser and Zannier. We proceed as follows:  from the above list $H_1,\ldots   H_t$ for each one we pick one of the multiplicative conditions defining each of them. These define a finite set $\{ J_1,\ldots ,J_t\}$ of codimension $1$ subgroups so that every anomalous subvariety is contained in a coset of one of them. Crucially, these cosets can be chosen defined over $L$.

\par
Let $W$ be anomalous as above. Then there is a condimension $1$ subgroup $J_j$ from the above finite list defined by a nonzero  $D(k_2+1)$ integer
vector,
$$\overline{ \alpha_j} = \alpha_{j1},\ldots ,\alpha_{jD(k_2+1)}$$
and $\theta_W\in L$ such that the following relation holds


\begin{equation}
\label{identity} \overline{w}^{\overline{\alpha_{j}}}=\theta_W
\end{equation}
for all $\overline{w} \in W$. Notice that the finitely many
vectors $\overline{\alpha_1}, \ldots ,\overline{\alpha_t}$
depend only on the variety $V^{k_2+1}.$ Fix an  order on the finite set of
$\overline{\alpha_{j}}$'s, $j=1,\ldots , t$. To any subset of $S$
of cardinality $k_2+1$, $E=\{ \eta_1, \ldots ,\eta_{k_2+1}\}$,
where $\eta_1< \ldots <\eta_{k_2+1}$ with respect to the fixed
order on $S$,  we associate the anomalous variety of  the tuple
$e^{\overline{b}\eta_1}, \ldots , e^{\overline{b}\eta_{k_2+1}}$.

We now define a coloring of the subsets of $S$ of cardinality
${k_2+1}$. Consider the function
$$ \Phi :[S]^{k_2+1}\rightarrow \{ \overline{\alpha_{1}},  \ldots ,\overline{\alpha_t} \}$$
that associates to any set in $[S]^{k_2+1}$ the tuple
$\overline{\alpha_{j}}$ for the minimum  $j$ such that the
anomalous variety corresponding to the subset is included in a
coset (defined over $L$) of $J_j$.

By Ramsey's  Theorem there is an infinite set $T\subseteq S$  and
a fixed $j_0$ such that $ \Phi$ takes the constant value
$\overline{\alpha_{j_0}}$ on the set of $k_2 + 1$ cardinality
subsets of $T$. Let $F\in [T]^{k_2+1} $ and order the elements of
$F$ as $\epsilon_1<\ldots < \epsilon_{k_2+1}$, where $<$ is the
order of $T$ inherited from $S$. We write the $D(k_2+1)$-tuple
$\overline{\alpha_{j_0}}$  as the concatenation of two parts
$\overline{\alpha_{j_0}}=\overline{\alpha_{j_0+}}$ $
\overline{\alpha_{j_0-}}$, where the {\it minus} part denotes the
last block of $D$ elements.


\noindent {\bf Case a)}
$\overline{\alpha_{j_0}}_-\not=\overline{0}.$ We fix the first
${k_2}$ elements of $F$, $\epsilon_1<\ldots < \epsilon_{k_2}$, and
we consider all elements $s$ of $T$ greater than $\epsilon_{k_2}$.
There are infinitely many such elements $s$, and if we append the
$D$-tuple $b_1s, \ldots, b_Ds$  to  ${\overline{b}\epsilon_1}
\ldots {\overline{b}\epsilon_{k_2}}$, we get an element of
$[T]^{k_2+1} $. We now exploit the indiscernibility of
$[T]^{k_2+1} $. For each $s$ as chosen, there is an element
$\theta_s $  in $L$ such that
\begin{equation}
\label{relation1} (e^{\overline{b}\epsilon_1} \ldots
e^{\overline{b}\epsilon_{k_2}})^{\overline{\alpha_{j_0}}_{+}}(e^{\overline{b}s}
)^{\overline{\alpha_{j_0}}_-} =\theta_s.
\end{equation}

Let $A_F= (e^{\overline{b}\epsilon_1} \ldots
e^{\overline{b}\epsilon_{k_2}})^{\overline{\alpha_{j_0}}_{+}}$.
Notice that the inner product $\overline{b}\cdot
\overline{\alpha_{j_0}}_-\not= 0$ since
$\overline{\alpha_{j_0}}_{-}\not= \overline{0}$ and the $\overline{b}$
are linearly independent over $\mathbb Q$. (Recall that we always assume $s$ is nonzero). So
$$(e^{\overline{b}s} )^{\overline{\alpha_{j_0}}_-}
=\frac{\theta_s}{A_F} \in L(A_{F}) $$ for each fixed $s$.

Then the
transcendence degree of $$\{ e^{(\overline{b}\cdot \overline{\alpha_{j_0}}_{-})s} : s\in T \backslash
 \{\epsilon_1, \ldots, \epsilon_{k_2}\}\}$$ over $\mathbb Q$ is bounded by the
transcendence degree of $L$ and $A_F$.  Appealing to Schanuel's conjecture we get
the finiteness of the linear dimension of the $\mathbb Q$-space
generated by  $\{ (\overline{b}\cdot \overline{\alpha_{j_0}}_{-})s$ : $s\in T \backslash
 \{\epsilon_1, \ldots, \epsilon_{k_2}\} \}$.  Clearly, then the set $ T \backslash
 \{\epsilon_1, \ldots, \epsilon_{k_2}\}$ is finite dimensional over $\mathbb Q$.

\par
\noindent {\bf Case b)} $ \overline{\alpha_{j_0}}_-=\overline{0}$. We shift to the
next block to the left in $ \overline{\alpha_{j_0}}$ not identically zero
suppose this corresponds to $\ell $, with $\ell \leq k_2$. As before we make the
corresponding  coordinate in the $\ell$th position in $\epsilon_1<\ldots < \epsilon_{\ell }$ vary over all
elements of $T$ strictly greater than $\epsilon_1<\ldots <
\epsilon_{\ell -1}$ and completing the $k_2+1$ tuple respecting the
order of $T$. We argue then as before. \qed
\\

An immediate consequence of the finite dimensionality of  $S$  is  the following corollary which can be viewed as a
multiplicative version of the statement of the Main Lemma.

\begin{corollario}\label{hull}
Let $\widehat{G}$ be the divisible hull of $G$, the group generated by all $e^{\mu_js}$'s where $s \in S$ and $j=1,\ldots ,N$.   Then $\widehat{G}$ has finite rank.\\
\end{corollario}
A basic result on linear functions on finite rank groups that will be relevant in the remaining part of the proof is due to Evertse, Schlickewei and Schmidt (see \cite{evertse}).

 We recall that a solution $(\alpha_1, \ldots, \alpha_n)$ of a linear equation
  \begin{equation}
 \label{linearomog}
 a_1x_1 + \ldots + a_nx_n = 1
 \end{equation}
  over a field $K$ is non degenerate if for every  proper non empty subset $I$ of $\{1,
\ldots, n\}$ we have  $\sum_{i \in I}a_i\alpha_i \not= 0$.

In our context we will be interested in solving the linear equation (\ref{linearomog}) in units of the field. Hence it is natural to consider equations of the form $$a_1x_1 + \ldots + a_nx_n  = 0$$ instead than (\ref{linearomog}).

\begin{lemma}[\cite{evertse}]
\label{viola} Let $K$ be a field of characteristic $0$, $n$ a
positive integer, and  $\Gamma$  a finitely generated subgroup of
rank $r$ of the multiplicative group $( K^{\times})^{n}.$ There
exists a positive integer $R = R(n, r)$ such that for any non zero
$a_1, \ldots, a_n$ elements in $ K,$ the equation
\begin{equation}\label{equazioneviola}
a_1x_1 + \ldots + a_nx_n  = 1 \end{equation} does not have more
than $R$ non degenerate solutions $(\alpha_1, \ldots, \alpha_n)$
in $\Gamma$.\\
\end{lemma}

We now apply this result to our context. Let $p \in \Bbb N,$ be
the linear dimension of $S,$ and  $\{ s_1,\ldots ,s_{p}\}$ be a $\mathbb Q$-basis of $S$.
For any $s \in S$ we have

\begin{equation}\label{soluzioni}
s= \sum_{l =1}^{p}c_{l}s_l
\end{equation}
where $c_{l} \in \Bbb Q.$ Substituting the expression of $s$
as in (\ref{soluzioni}) in $f$ we have

\begin{equation}\label{polexp}
0 = f(s) = \lambda_1 e^{\mu_1(\sum_{l =1}^{p}c_{l}s_l)} +
\ldots + \lambda_N e^{\mu_N(\sum_{l =1}^{p}c_{l}s_l)} =
\sum_{j=1}^{N} \lambda_j \prod_{l=1}^{p}
(e^{\mu_js_l})^{c_{l}}
\end{equation}

Any solution $s \in S$  produces a solution
$\overline \omega$  of the linear equation associated to $f$,
\begin{equation}\label{lineare}
\lambda_1 X_1 + \ldots + \lambda_N X_N = 0 \end{equation} where
$\omega_i = e^{\mu_i(\sum_{l =1}^{p}c_{l}s_l)},$ $i = 1, \ldots, N$
and  $\overline \omega \in \widehat{G}$ (a subgroup of $\Bbb (\Bbb C^{*})^N$, see Corollary \ref{hull}).

Since the coefficients of $f$ are nonzero, we can transform this equation to the form of the unit equation by
 replacing $\lambda_N$ by $-\lambda_N$, and multiplying throughout by $(-\lambda_N)^{-1}e^{-\mu_Ns}.$

\begin{lemma}\label{soluzdistinte} Suppose $f(z) = \lambda_1e^{\mu_1 z} + \ldots + \lambda_Ne^{\mu_N z}$
is not simple, and $s_{1}, s_{2}$ are two distinct solutions of
$f.$ Then the solutions of (\ref{lineare}) generated by $s_{1}$
and $s_{2}$ are different.\end{lemma}

\noindent \emph{Proof:} Let $\overline \omega = \omega_1, \ldots, \omega_N$
and $\overline \xi = \xi_1, \ldots, \xi_N$ be the
solutions of (\ref{lineare}) corresponding to $s_{1}$ and $s_{2}$,
respectively. If
$$(\omega_1, \ldots, \omega_N) = (\xi_1, \ldots, \xi_N),$$ then for $j =
1, \ldots, N,$
$$\prod_{l=1}^{p} (e^{\mu_js_l})^{c_{1,l}} = \prod_{l=1}^{p} (e^{\mu_js_l})^{c_{2,l}},$$
iff
$$\prod_{l=1}^{p} (e^{\mu_js_l})^{c_{1,l} - c_{2, l}} = 1$$
iff
$$\mu_j \sum_{l=1}^{p} s_l(c_{1,l} - c_{2, l}) \in 2\pi i \Bbb Z.$$
So, for any $j=1, \ldots, N$ we have
$$\sum_{l=1}^{p} s_l(c_{1,l} - c_{2, l}) = \frac{2\pi i}{\mu_j}h_j$$
where $h_j \in \Bbb Z.$ This implies $$\frac{2\pi i}{\mu_1} h_1 =
\frac{2\pi i}{\mu_2} h_2 = \ldots = \frac{2\pi i}{\mu_N} h_N.
$$
So we can write any exponents $\mu_j$ in the polynomial $f(z)$ in
terms of $\mu_1,$ i.e.
$$\mu_2 = \frac{\mu_1}{h_1} h_2$$
$$\mu_3 = \frac{\mu_1}{h_1} h_3$$
$$\ldots$$
$$\mu_N = \frac{\mu_1}{h_1} h_N.$$
If $\alpha = \frac{\mu_1}{h_1}$ then $f(z)$ is a polynomial in
$e^{\alpha z},$ i.e. $f$ is a simple polynomial. We get a contradiction since $f$ is not  simple.\\


We restate the remaining case of Shapiro's Conjecture.

\begin{thm}\label{shapiro} (SC) Let $f(z)$ be an irreducible polynomial and suppose the following system
\begin{equation}\label{sistema1}
\left\{ \begin{array}{lll} f(z) = \lambda_1e^{\mu_1 z} + \ldots + \lambda_Ne^{\mu_N z} & = 0 &\\
 g(z) = l_1e^{m_1 z} + \ldots + l_Me^{m_M z} & = 0 &
\end{array} \right.\\
\end{equation}
has infinitely common zeros. Then $f$ divides $g$.
\end{thm}

\noindent {\it Proof:}  We will use induction on the length of the
polynomial $g(z).$ Without loss of generality we may assume $N, M > 2,$ and $g$ not
simple otherwise we
would be in Case 1 solved by van der Poorten and Tijdeman.\\
Consider the linear equation  associated to $g(z)=0$,
\begin{equation}\label{g}
l_1X_1 + \ldots + l_MX_M = 0.
\end{equation}

We can transform this equation to a unit equation as in Lemma \ref{viola}. Lemma \ref{soluzdistinte}
implies that equation (\ref{g}) has infinitely many solutions  $\overline \alpha = (\alpha_1, \ldots, \alpha_M),$ where
$\alpha_t = e^{m_t(\sum_{\ell =1}^{p}c_{\ell}s_{\ell})}$ (each one
generated by $s$, a solution of (\ref{sistema1})). From Lemma
\ref{viola} it follows that all but finitely many of them are degenerate.
\par By the Pigeonhole Principle there exists a proper subset $I = \{i_1, \ldots, i_r\}$  of $\{1,
\ldots, M\}$ such that
\begin{equation}\label{sist6}
l_{i_1}X_{i_1} + \ldots + l_{i_r}X_{i_r}= 0
\end{equation}
has infinitely many zeros  of the right form. Notice that $I$ has at least three elements since we are assuming that $g$ is not a simple polynomial. \\

It is useful to write $g(z) = g_1(z) + g_2(z)$, where $g_1(z)= l_{i_1}
e^{m_{i_1}z}+ \ldots + l_{i_r}e^{m_{i_r}z},$ and  $g_2(z)=g(z) - g_1(z)$. The polynomial
$g_1$ has infinitely
many common zeros with $f(z).$ Also, $g_2(z)$ has infinitely many common zero with $f(z)$. Both $g_1(z)$ and $g_2(z)$ have  lengths strictly less than $M.$

By inductive hypothesis and by the
irreducibility of $f$, we have that $f$ divides $g_1$ and $f$
divides $g_2,$ and hence $f$ divides $g.$ So the proof  is
completed.\\
\qed






\begin{thebibliography}{99}


\bibitem{bombieri} E. Bombieri, D. Masser and U. Zannier: \emph{Anomalous Subvarieties Structure Theorems and Applications},
International Mathematics Research Notices 2007, (2007), 1-33 .



\bibitem{paolagiusy} P. D'Aquino and G. Terzo: \emph{A theorem of complete reducibility  for exponential polynomials}, submitted, (2011).


\bibitem{vanderpoorten} A. J. van der Poorten: \emph{Factorisation in fractional powers},
Acta Arithmetica 70, (3),  (1995), 287-293.

\bibitem{vdpoorteneverest} A. J. van der Poorten and G. R. Everest:
\emph{Factorisation in the ring of exponential polynomials},
Proceedings of the American Mathematical Society, 125, (5),
(1997), 1293-1298.
\bibitem{vanderpoortentijdeman} A. J. van der Poorten and R.
Tijdeman:\emph{On common zeros of exponential polynomials},
L'Enseignement Math\'{e}matique, 21, (1975), 57-67.

\bibitem{evertse} J. H. Evertse, H. P. Schilickewei and W. M.
Schmidt, \emph{Linear equations in variables which lie in a
multiplicative group,} Annals of Mathematics, 155 (2002), 807-836.

\bibitem{gourin} E. Gourin: \emph{On irreducible polynomials in several variables
which become reducible when the variables are replaced by powers
of themselves}, Transactions of the American Mathematical Society
32, (1930), 485-501.


\bibitem{KMO}
J. Kirby, A. Macintyre, and A. Onshuus: \emph{The algebraic numbers definable
in various exponential fields}, to appear in J. Inst. Math. Jussieu, June
2010

\bibitem{macoll} L. A. MacColl: \emph{A factorization theory for polynomials in $x$ and in functions
 $e^{\alpha x}$}, Bulletin of the American Mathematical Society, (41), (1935), 104-109.



\bibitem{angusfree} A. Macintyre: \emph{Exponential Algebra}, in Logic and Algebra. Proceedings of the
international conference dedicated to the memory of Roberto
Magari, (A. Ursini et al. eds), Lecture Notes in Pure Applied
Mathematics 180, (1991), 191-210.

\bibitem{MacWilkie} A. Macintyre and A. Wilkie: \emph{On the decidability of the real exponential field},
in Kreiseliana: about and around Georg Kreisel, A. K. Peters
(1996), 441-467.


\bibitem{Marker} D. Marker: \emph{A remark on Zilber's pseudoexponentiation}, The Journal of Symbolic
Logic, 71, (3), (2006), 791-798.

\bibitem{ritt} J.F. Ritt: \emph{A factorization theorem of functions $\sum_{i=1}^na_ie^{\alpha_iz}$}, Transactions of the
American Mathematical Society 29, (1927), 584-596.

\bibitem{ritt1} J.F. Ritt: \emph{On the zeros of exponential
polynomials}, Transactions of the American Mathematical Society, 31,
(1929), 680-686.


\bibitem{shapiro conj} H. S. Shapiro: \emph{The Expansion of mean-periodic functions in series of
exponentials}, Comm. Pure and Appl. Math. 11, (1958), 1-21.



\bibitem{GiusyEuler} G. Terzo: \emph{Some Consequences of Schanuel's Conjecture in
exponential rings} Communications in Algebra 36, (3), (2008),
1171-1189.

\bibitem{Zilber} B. Zilber: \emph{Pseudo-exponentiation on
algebraically closed fields of characteristic zero}, Annals of
Pure and Applied Logic, 132, (1), (2005), 67-95.


\end{thebibliography}
\end{document}